\documentclass[11pt,a4paper]{article}
\usepackage{graphicx}
\usepackage{amsmath, amssymb}

\usepackage[T1, OT1]{fontenc}

\DeclareTextCommand{\Dj}{OT1}{%
\raisebox{-0.1ex}{\scalebox{0.75}[1.4]{--}}\kern-.4em D%
}
\DeclareTextCommand{\dj}{OT1}{%
\raisebox{0.4ex}{\scalebox{0.75}[1.6]{--}}\kern-.5em d%
}

  \newtheorem{theorem}{Theorem}[section]
  \newtheorem{corollary}{Corollary}[theorem]
  \newtheorem{lemma}[theorem]{Lemma}
  \newtheorem{proposition}[theorem]{Proposition}
  \newtheorem{definition}[theorem]{Definition}

\providecommand{\keywords}[1]{\textbf{\textit{Keywords:\,\,}} #1}
  \newcommand{\R}{{\cal R}}
  \newcommand{\Z}{{\cal Z}_0}
  
  \newcommand{\N}{{\cal N}}
  
  \newcommand{\lL}{{\cal L}}
  
  \newcommand{\stp}{{\rm st}}
  \newcommand{\Rf}{{\R_{\rm fin}^*}}

\begin{document}

  \title{A nonstandard approach to Karamata uniform convergence theorem }
  \author{\v Zarko Mijajlovi\'c \\
          Faculty of Mathematics, University of Belgrade, \\
          Studentski trg 16,  11000 Belgrade, Serbia \\
          zarkom@matf.bg.ac.rs\\  \\
          Danijela Brankovi\'c\\
          School of Electrical Engineering, University of Belgrade, \\
          Bulevar kralja Aleksandra 73,  11000 Belgrade, Serbia \\
          danijela@etf.bg.ac.rs
         }

  \date{}
  \maketitle

  \begin{abstract}
    A nonstandard proof of a generalization of Karamata uniform convergence theorem for slowly varying
    functions is presented. Properties of a related operator $\mathcal{L}$ and its connection with slowly varying functions are discussed. 
  \end{abstract}
\noindent
 
\section{Introduction}
\label{intro}

  This work is inspired by the possibility of application of the theory of
  regularly varying functions in study of asymptotics of cosmological parameters.
  Our particular aim was to consider the asymptotics of the expansion scale factor $a(t)$ in
  the $\Lambda$CDM (Lambda cold dark matter) model, see \cite{liddle}. It appeared that the mathematical behavior of the
  scale factor $a(t)$ in certain epochs  of evolution of the Universe is very connected to the
  properties of regularly varying functions.
  Some works in this area are \cite{stern},  \cite{mijajlo2012}, \cite{mijajlo2015} and \cite{mijajlo2019}.
  For that cause, we extracted and studied certain general properties of regularly varying functions appearing in
  our cosmological studies  which might be of an interest by themselves.
  Obtained results on these properties are presented in this paper.

  We shall occasionally  use here the methods of nonstandard analysis, as we did in \cite{mijajlo2007} and \cite{mijajlo2007A}.
  Somewhat extended explanation of
  notions from nonstandard analysis is given so that a non-specialists in this area can read the paper, too.
  Only basic notions from model theory of first order logic will be assumed, see \cite{Keisler}. For
  more details of this subject one may consult \cite{Stroyan}, \cite{Davis} and \cite{mijajlo2014}.

  If $\R$ stands for the field of real numbers with some added functions and relations, then $\R^*$ denotes a
  nonstandard extension of $\R$. We remind that $\R \prec \R^*$, i.e. $\R$ is elementary embedded in $\R^*$.
  In other words all first order properties expressed in the expanded language
  $L= \{+,\cdot, \leq\} \cup \{\underline{s}\colon s\in S\}$, where $S$ is the set of
  constants from $\R$ and added functions and relations, are preserved from $\R$ to $\R^*$ and vice versa.
  We call it the transfer principle.
  Here $\underline{s}$ is the name of an entity $s$ from $S$. The symbol $\underline{s}$ is interpreted as $s$ in $\R$, i.e.
  $\underline{s}^\R= s$, while $s^*=\underline{s}^{\R^*}$. Due to  $\R \prec \R^*$
  if $\varphi$ is a sentence of $L$, then $\R\models \varphi$ if and only if $\R^* \models \varphi$, where $\models$ is
  the satisfaction relation.
  If $\varphi= \varphi(\underline{s}_1,\underline{s}_2,\ldots \underline{s}_n)$, we shall often write
  $\R^*\models \varphi(s_1^*,s_2^*,\ldots,s_n^*)$ instead of $\R^*\models \varphi$. For easier reading we shall omit
  in some cases the star in $s^*$ if this does not lead to ambiguity, for example in $+^*$, $\leq^*$, etc.

  We use the same symbols for the structures and  their domains,
  e.g. $\N$ denotes the set of natural numbers (non-negative integers) and the structure of natural numbers, too.
  Infinitesimal is an  element $\varepsilon\in \R^*$ that is infinitely close to $0$, i.e. for all $n\in \N^+$,
  $|\varepsilon| < 1/n$. The symbol $\mu(0)$ stands for the set of all
  infinitesimals   and is called the monad of zero.
  If $a-b$ is an infinitesimal, then we write $a\approx b$.
  An element $H\in \N^*\setminus \N$ is a positive infinite integer if for all $n\in \N$, $H>n$, while  $a\in \R^*$ is  finite
  if for all infinite positive integers $H$, $|a|<H$.
  Elements of $\R^*$ that are not finite are called infinite.
  Let $\Rf$ denote all finite elements of $\R^*$.
  Then for   $b\in\Rf$ there is $a\in \R$ and an infinitesimal $\varepsilon$ such that
  $b= a+ \varepsilon$. Then $\stp(b)=a$, where $\stp$ is the standard part
  function $\stp\colon \Rf\to \R$.
  We remind that $\stp$ is a homomorphism from the field
  $\Rf$ to $\R$ and for continuous functions as well.
  It is convenient to extend $\stp$ to infinite elements of $\R^*$, taking for
  positive infinite $H\in \R^*$, $\stp(H)= \infty$ and $\stp(-H)= -\infty$.
  Monad of $b\in \R$ is $\mu(b)= \mu(0)+b$.
  It is convenient to use the quasi-order $\lesssim$ on $\R^*$ defined by
  $x\lesssim y$, $x,y\in \R^*$ if and only if $x\leq^* y$ or $x\approx y$.
  If $\varphi(x)$ is a predicate formula which defines a set $X\subseteq \R$,
  then $\varphi^*(x)$, obtained by starring entities over $\R$ appearing in $\varphi$,
  defines an internal set $X^*\subseteq \R^*$ associated to $X$.
  On internal subsets can be defined a finite-additive measure which can be extended by use of Caratheodory extension theorem to
  $\sigma$-additive measure on $\R^*$. This measure is called Loeb measure which is naturally related to
  Lebesgue measure.
  For other notation and terminology  see \cite{Keisler}, \cite{Stroyan}, \cite{mijajlo2014}.

  The symbol  $\bigwedge$ denotes the universal quantifier, while $\bigvee$
  stands for the  existential  quantifier. For example, for
  $F\colon \R\times I \to \R$, $I\subseteq \cal R$,  we have

  \begin{multline}
    \lim_{x\to\infty} F(x,u)= 0\quad  {\rm for\,\, all}\quad u\in I\quad  {\rm if\,\, and\,\, only\,\, if}\quad \\
    \R\models \bigwedge_{u\in I}\bigwedge_{\varepsilon>0}\bigvee_{x_0} \bigwedge_{x>x_0} |F(x,u)|<\varepsilon.
  \end{multline}

  \begin{multline}\label{unicov}
    \lim_{x\to\infty} F(x,u)= 0\quad  {\rm uniformly\,\, for}\quad u\in I\quad  {\rm if\,\, and\,\, only\,\, if}\quad \\
    \R\models \bigwedge_{\varepsilon>0}\bigvee_{x_0} \bigwedge_{x>x_0}\bigwedge_{u\in I} |F(x,u)|<\varepsilon.
  \end{multline}

  A real function $F$ is said to be regularly varying at infinity if it is real-valued, positive and
  measurable on $[a,\infty)$, for some $a > 0$, and if for each $\lambda > 0$
  \begin{equation} \label{slovar}
    \lim_{x\to\infty}\frac{F(\lambda x)}{F(x)} = \lambda^\rho
  \end{equation}
  for some $\rho$,  $-\infty < \rho < \infty$.

  Number $\rho$ is called the index of regular variation. If $\rho=0$, $F$ is called slowly varying, or SV  function.
  Notions of slowly varying functions and regular variations were introduced by Jovan Karamata \cite{karamata}.
  Books \cite{bingham} and \cite{Seneta}
  give detailed exposition of the theory of regular variation and slowly varying functions.
  The following theorem which refers to slowly varying functions is fundamental in this theory.

  \begin{theorem}\label {karamata}
  {\rm (The Uniform Convergence Theorem, J. Karamata)}. If $F$  is a slowly varying function,
   then for every fixed $[a,b]$, $0 < a < b < \infty$, the relation (\ref{slovar}) holds uniformly with respect to $[a,b]$.
  \end{theorem}

  We shall prove a generalization of the linear variant of this theorem using nonstandard methods.
  The linear form is obtained by  transformation $f(x) = \ln F(e^x)$.

\section{Regular variation in nonstandard analysis}
\label{regular}

  First we prove a lemma on uniform convergence relative to a set $I$.
    \medskip

  \begin{lemma} \label{unic}
     Let $F\colon \R\times I \to \R$,   $I\subseteq \R$. Then

     $\lim_{x\to \infty} F(x,u)= 0$ uniformly with respect to $u\in I$ if and only if

     for all positive infinite $x\in \R^*$, $u\in I^*$, $F^*(x,u) \approx 0$.

  \end{lemma}

  \noindent{\bf Proof}\quad
    ($\Rightarrow$) Suppose $\lim_{x\to \infty} F(x,u)= 0$ uniformly for $u\in I$, i.e.
    (\ref{unicov}) holds. Let $\varepsilon\in \R^+$ be arbitrary and $x_0\in \R$ so that
    $\R\models \bigwedge_{x>x_0}\bigwedge_{u\in I} |F(x,u)|<\varepsilon$.
    By transfer principle,
     $\R^*\models \bigwedge_{x>x_0}\bigwedge_{u\in I^*} |F^*(x,u)|<\varepsilon$.
    Hence, for all positive infinite $x\in \R^*$,
     $\R^*\models \bigwedge_{u\in I^*} |F^*(x,u)|<\varepsilon$, and so for all
     $u\in I^*$, $|F^*(x,u)|<\varepsilon$. As $\varepsilon\in \R^+$ was chosen  arbitrarily,
     it follows $F^*(x,u)\approx 0$.

     \noindent ($\Leftarrow$)  Suppose for all positive infinite $x\in \R^*$, $u\in I^*$, $F^*(x,u) \approx 0$.
     Then for   arbitrary  $\varepsilon\in \R^+$,
     all positive infinite $x$ and all $u\in I^*$, $|F^*(x,u)|<\varepsilon$.
     Choose any positive  infinite $x_0$. As $x> x_0$ is also positive infinite it follows
     $\R^*\models \bigwedge_{x>x_0}\bigwedge_{u\in I^*} |F^*(x,u)|<\varepsilon$.
     Hence for all $\varepsilon\in \R^+$,
     $\R^*\models \bigvee_{x_0}\bigwedge_{x>x_0}\bigwedge_{u\in I^*} |F^*(x,u)|<\varepsilon$, i.e.
     $\R\models \bigvee_{x_0}\bigwedge_{x>x_0}\bigwedge_{u\in I} |F(x,u)|<\varepsilon$.
     It follows
     $\R\models \bigwedge_{\varepsilon>0} \bigvee_{x_0}\bigwedge_{x>x_0}\bigwedge_{u\in I} |F(x,u)|<\varepsilon$,
     so $\lim_{x\to \infty} F(x,u)= 0$ uniformly in respect to $u\in I$.
     \hfill $\Box$

     \begin{corollary} \label{cor} $\lim_{x\to\infty} F(x,u)= 0$ is not uniformly convergent relative to $u\in I$
     if and only if there is a positive infinite $x^*$, $u^*\in I^*$ and $a\in \R^+$ such that $|F^*(x^*,u^*)|>a$.
     \end{corollary}

     In contrast to this corollary, note that for the ordinary convergence, $\lim_{x\to\infty} F(x,u)= 0$
     if and only if
     for all positive infinite $x\in \R^*$, $u\in I$, $F^*(x,u) \approx 0$.
     \smallskip

     Now we  prove a generalization of Karamata uniform convergence theorem for slowly varying
     functions. In the next proof we assume the notion and properties of Loeb measure, a natural extension
     of Lebesgue measure into the nonstandard universe.
     If $A\subseteq\R$ is a measurable, $\lambda(A)$ denotes Lebesgue measure of $A$,
     while $\ell(A)$ denotes  measure of a Loeb measurable set $A\subseteq\R^*$  in the nonstandard universe.
     The following well-known Fisher's lifting theorem gives a connection between Lebesgue and Leobe measure
     and basically states that $\stp^{-1}\colon \R \to \Rf$ preserves measure.
     \begin{theorem}
       If $A$ is a Lebesgue measurable subset of a finite closed interval of  $\R$, then
       \begin{equation}
          \lambda(A)=  \ell(\stp^{-1}(A)).
       \end{equation}
     \end{theorem}

     Now we proceed to the proof of the main theorem in Sec.~\ref{regular}, a generalization of
     Karamata uniform convergence theorem.

     \begin{theorem} Assume  $H\colon S\times \R \to \R^+$ is a measurable,
     where $S= [d,\infty]_\R$ for some $d\in \R$, $I=[0,1]_\R$ and suppose $H$ satisfy the following inequality on its domain:
     \begin{equation}\label{HI}
          H(x,u)\leq H(x+u,v) + H(x,u+v).
     \end{equation}
     Further, let $m\colon S\to \R^+$ be a measurable and nondecreasing function and $G(x,u)= H(x,u)m(x)$.
     Then, if $\lim_{x\to \infty} G(x,u)= 0$ for all $u\in \R$, then this convergence is uniform relative to $u\in I$.
     \end{theorem}

     \noindent{\bf Proof} 
     Assume $\lim_{x\to \infty} G(x,u)= 0$ for all $u\in \R$, but this convergence is
     not uniform relative to $u\in I$. Then by Corollary \ref{cor} there are positive infinite $x_0^*\in \R^*$, $u_0^*\in I^*$
     and $a\in \R^+$ such that $G^*(x_0^*,u_0^*)> a$.
     Let $[0,2]^*= [0,2]_{\R^*}$ and define
     \begin{equation}
       \begin{split}
         {U} &= \{u\in [0,2]^*\, |\, G^*(x_0^*,u) < a/3 \}  \\
         {V} &= \{v\in [0,2]^*\, |\, G^*(x_0^* + u_0^*,v) < a/3 \}.
       \end{split}
     \end{equation}
     As sets $U$ and $V$ are internal and $G$ is $\lambda$-measurable, $U$ and $V$
     are $\ell$-measurable. For positive infinite $x$ and $u\in [0,2]$ we have $G^*(x,u)\approx 0$, hence
     $[0,2]\subseteq U\cap V \subseteq [0,2]^*$. Then  Fisher's theorem implies
     \begin{equation}
       \ell(U) =2,\quad  \ell(V)=2.
     \end{equation}
     Let $V_0= V + u_0^*$. As  measure $\ell$ is invariant under translation,
     we have $\ell(V_0)= \ell(V)$, i.e.
     \begin{equation}
         \ell(V_0)= 2.
     \end{equation}
     Further, $U, V_0 \subseteq [0,3]^*$, hence $U\cap V_0\not=\emptyset$,
     as otherwise
     \smallskip

     $3=\ell([0,3]) \geq \ell(U \cup V_0)= 4$, a contradiction.
     \smallskip

     So, let $b^*\in U\cap V_0$. Hence, $b^*\in U$ and for some $v^*\in V$,
     $b^* = v^* + u_0^*$ and so $v^*= b^* - u_0^*$. Then
     \smallskip

     As $v^*\in V$, we have $G^*(x_0^* + u_0^*,v^*) < a/3$,

     As $b^* \in U$, we have $G^*(x_0^*,b^*)<a/3$.

     Hence,
     \begin{equation}
        H^*(x_0^* + u_0^*,v^*)m(x_0^* + u_0^*) < a/3
     \end{equation} and
     \begin{equation}
        H^*(x_0^*,b^*)m(x_0^*) < a/3.
     \end{equation}
     By the inequality (\ref{HI}) we have
     \begin{equation}
       \begin{split}
             H^*(x_0^*,u_0^*) &\leq H^*(x_0^*+u_0^*,v^*) + H^*(x_0^*,u_0^* + v^*)  \\
                              &= H^*(x_0^*+u_0^*,v^*) + H^*(x_0^*,b^*)
       \end{split}
     \end{equation}
      As $m(x)$ is nondecreasing we have $ m(x_0^*) \leq m(x_0^*+u_0^*)$, so
     \begin{equation}
          H^*(x_0^*,u_0^*)m(x_0^*) \leq H^*(x_0^* + u_0^*,v^*)m(x_0^*+u_0^*) + H^*(x_0^*,b^*)m(x_0^*),
     \end{equation}
     i.e.
     \begin{equation}
       a< G^*(x_0^*,u_0^*) \leq G^*(x_0^* +u_0^*,v^*) + G^*(x_0^*,b^*) < 2a/3,
     \end{equation}
     a contradiction.
     \hfill $\Box$

     The uniform convergence is preserved under translation and homothety, hence we have:
     \begin{corollary}
       The previous theorem is still true if  the interval $I=[0,1]_\R$ is replaced by any finite
       interval $I'= [a,b]_\R$, $0<a<b$.
     \end{corollary}

     If we take $H(x,u) = |f(x+u) - f(x)|$ and $m(x)=1$, where $f(x)= \ln(F(e^x))$ and $F(x)$  is a slowly varying function,
     we immediately obtain Karamata theorem \ref{karamata}.


\section{Operator $\lL$}
\label{operator}

  We introduce  operator $\lL$,  which have an important role in the analysis of regular variation.
  The  operator $\lL$ may be defined on the set of Lebesgue integrable functions,
  but due to the nature of physical  parameters  that are studied in this paper,
  our attention will be turned only towards to its restriction  to at least twice differentiable
  real functions, i.e to the space $C^2(R)$.
  \begin{definition}\label{op}\,\,
    ${\lL}(h)(x) = \displaystyle{\frac{1}{\ln(x)}\int_{x_0}^x \, \frac{h(t)}{t}dt}$,
    \quad $x_0>1$, $h\in C^2(R)$.
  \end{definition}

  As we are interested in asymptotics at infinity, the exact value of $x_0$ is not of some importance.
  We may even assume,  with a proper adaptation of the argument function $h$,  that $x_0=1$.
  Namely, if the function $h(x)$ in the above definition of $\lL$ is bounded in some neighborhood of $1$,
  what is in this paper almost of the only interest, then  by the l'Hopital's rule
  \begin{equation}
     \lim_{x\to 1} {\lL}(h)(x)= h(1),
  \end{equation}
  so we can take in the above definition  $x_0=1$. From now on we assume $x_0=1$, if it is not otherwise specified.
  Obviously, $\lL$ is a linear operator over the space $C^2(R)$. This operator  has many interesting properties
  and some of them reflect more or less well-known
  theorems on regularly varying functions.

  It is convenient to denote by $\cal R_\alpha$ the class of regularly varying functions of index $\alpha$.
  Hence ${\cal R}_0=$ SV is the class of all slowly varying functions.
  By ${\cal Z}_0$ we shall denote the class of zero functions at $\infty$, i.e.
  $\varepsilon\in {\cal Z}_0$ if and only if $\displaystyle\lim_{t\to \infty} \varepsilon(t)=0$.

  \begin{theorem} \label{thmL}
  Let $\lL'$ denote the restriction of $\lL$ to the appropriate domain. Then
  \begin{enumerate}
    \item $\lL'\colon \Z\to \Z$.
    \item $\lL'\colon \R_0\to \R_0$.
    \item $\lL'\colon \R_{\alpha}\to \R_{\alpha}$,\quad $\alpha\in R$.
    \item $\lL'\colon B(R)\to B(R)$, where $B(R)$ is the set of real bounded functions.
  \end{enumerate}

  \end{theorem}

  \noindent{\bf Proof}.
   The statement 1. may be obtained by use of  l'Hopital's rule.
  Statement 2.  immediately follows from the statements 1.5.9a and 1.5.9b in \cite{bingham}.
  Statement 3. follows from theorems 1.5.10 and 1.5.11 in \cite{bingham}.
  Finally, if $h$ is bounded, then there are $a$ and $b$ such that $a\leq h(x) \leq b$,
  $x > 1$. Hence $a\leq \lL(h)(x) \leq b$, for $x > 1$.
  \hfill $\Box$
  \medskip

  The following proposition gives us one interesting property of the linear operator $\mathcal{L}$.

  \begin{proposition}
   Linear operator $\mathcal{L}$ is invertible.
  \end{proposition}

  \noindent{\bf Proof}.
  Let $f(x)=\mathcal{L}(g)(x)$, where $g\in C^2(R)$. By definition of the operator $\mathcal{L}$, it follows
  \begin{equation}\label{df}
  f(x)\ln(x) = \displaystyle{\int_1^x \, \frac{g(t)}{t}dt}.
  \end{equation}
  Differentiating the equation (\ref{df}) over variable $x$, we infer
  \begin{equation}\label{dopuna}
  \frac{f(x)}{x}+\dot{f}(x)\ln(x)=\frac{g(x)}{x}.
  \end{equation}
  From the equation (\ref{dopuna}) we obtain
  \begin{equation}
   g(x)=f(x)+x\dot{f}(x)\ln(x),
  \end{equation}
  what proves the proposition.
  \hfill $\Box$
  \medskip

  Limit  $\lim_{x\to\infty}\lL(h)(x)$ acts as a proper generalization of
  $\lim_{x\to\infty} h(x)$. Namely,
  by easy application of l'Hopital's rule we have
  \begin{proposition}\label{Llim} Let $c$ be a real number and
     suppose $\lim_{x\to\infty}h(x)= c$. Then $\lim_{x\to\infty}\lL(h)(x)= c$.
  \end{proposition}
  The  example $h(x)= \sin(x)$ shows that $\lim_{x\to\infty}\lL(h)(x)$ is a proper extension of ordinary limit.
  Namely, $h(x)$ diverges at infinity, while $\lim_{x\to\infty}\lL(h)(x) $ converges.

  If a slowly varying function $L(x)$ is given by integral representation (see \cite{karamata})
  \begin{equation*}
      L(x)= g(x)\displaystyle e^{\int_{x_0}^x \frac{\varepsilon(t)}{t}dt},
           \quad {\rm where}  \hspace{2mm} \varepsilon\in \Z \hspace{2mm} {\rm and}\hspace{2mm} g(x)\to g_0\,\,\,
      {\rm as}\,\, x\to \infty,
  \end{equation*}
  then
  \begin{equation*}
      \frac{\ln(L(x))}{\ln(x)} = \frac{\ln(g(x))}{\ln(x)} + \lL(\varepsilon)(x),
  \end{equation*}
  wherefrom by previous two propositions
    \begin{equation}
      \frac{\ln(L(x))}{\ln(x)}\to 0 \,\,\,  {\rm as}\,\, x\to \infty.
  \end{equation}
  This statement can be  found already in \cite{Seneta}.
  Now we give an application to asymptotics  of certain integrals.


  \begin{theorem}
    Let $h(x)$ be a positive and $M$-bounded Lebesgue integrable function, $M>0$ and $\lambda >1$.
    Assume
    \begin{equation}\label{lcond}
      \lim_{x\to\infty} (\lL(h)(\lambda x) - \lL(h)(x))\ln(x)=0
    \end{equation}
    uniformly with respect to $\lambda$. Then
    \begin{equation}
      \int_x^{\lambda x} \frac{h(t)}{t}dt \approx \frac{\ln(\lambda)}{\ln(\lambda) + \ln(x)}\int_1^x \frac{h(t)}{t}dt.
    \end{equation}
  \end{theorem}

 \noindent{\bf Proof} 
  Let
    \begin{equation}
      I= \frac{\ln(x)}{\ln(\lambda) + \ln(x)}\int_1^{\lambda x} \frac{h(t)}{t}dt - \int_1^x \frac{h(t)}{t}dt.
    \end{equation}
    First observe that
  \begin{equation}
    \begin{split}
      \frac{\ln(x)}{\ln(\lambda) + \ln(x)}&\int_1^{\lambda x} \frac{h(t)}{t}dt= \\
                                          &\int_1^x \frac{h(t)}{t}dt + \int_x^{\lambda x} \frac{h(t)}{t}dt -
                                           \frac{\ln(\lambda)}{\ln(\lambda) + \ln(x)}\int_1^{\lambda x} \frac{h(t)}{t}dt,
    \end{split}
  \end{equation}
  so
  \begin{equation}\label{int}
    I=\int_x^{\lambda x} \frac{h(t)}{t}dt - \frac{\ln(\lambda)}{\ln(\lambda) + \ln(x)}\int_1^{\lambda x} \frac{h(t)}{t}dt.
  \end{equation}
  Further,
  \begin{equation}
     \int_1^{\lambda x} \frac{h(t)}{t}dt= \int_1^x \frac{h(t)}{t}dt + \int_x^{\lambda x} \frac{h(t)}{t}dt
  \end{equation}
  and
  \begin{equation}
     \int_x^{\lambda x} \frac{h(t)}{t}dt \leq M\int_x^{\lambda x} \frac{dt}{t} = M\ln(\lambda).
  \end{equation}
  Hence,
  \begin{equation}\label{fra}
    \frac{\ln(\lambda)}{\ln(\lambda) + \ln(x)}\int_x^{\lambda x} \frac{h(t)}{t}dt \leq
        \frac{M\ln(\lambda)^2}{\ln(\lambda)+\ln(x)} \to 0 \quad {\rm as}\quad  x\to\infty.
  \end{equation}

  Hence, by (\ref{int}) and (\ref{fra})
  \begin{equation}
    I= \int_x^{\lambda x} \frac{h(t)}{t}dt - \frac{\ln(\lambda)}{\ln(\lambda) + \ln(x)}\int_1^x \frac{h(t)}{t}dt
      + \varepsilon(x),
  \end{equation}
  where $\varepsilon(x) \to 0$ as $x\to\infty$. By (\ref{lcond}) and the assumption that this convergence is
  uniform with respect to $\lambda$, there is $\xi(x)$ which does not depend on $\lambda$ so that $I=\xi(x)$ and $\xi(x) \to 0$ as $x \to \infty$.
Hence
  \begin{equation}
     \int_x^{\lambda x} \frac{h(t)}{t}dt -
        \frac{\ln(\lambda)}{\ln(\lambda) + \ln(x)}\int_1^{x} \frac{h(t)}{t}dt +\varepsilon(x) = \xi(x),
  \end{equation}
  so
  \begin{equation}
      \int_x^{\lambda x} \frac{h(t)}{t}dt \approx \frac{\ln(\lambda)}{\ln(\lambda) + \ln(x)}\int_1^x \frac{h(t)}{t}dt.
  \end{equation}
   \hfill $\Box$

   The following theorem gives us a better insight in the convergence of the limit
   appearing in the previous theorem.

   \begin{theorem} \label{thmlm}
     Let $f(x)$ be a measurable function defined on $[a,\infty)$ for some real number $a$
     and $S\subseteq R^+$ a measurable set of positive measure. If for all $\lambda\in S$
     \begin{equation}\label{llim}
       \lim_{x\to\infty} (f(\lambda x) - f(x))\ln(x)=0,
     \end{equation}
     then {\rm (\ref{llim})} holds for all $\lambda\in R^+$.
   \end{theorem}
   In the proof of the theorem we follow ideas presented in \cite{Seneta} and it is
   achieved by proving next lemmas.

   \begin{lemma} \label{lemmaA}
     Let $f(x)$  and $S$ be as in the theorem \ref{thmlm} and suppose {\rm (\ref{llim})} for
     all $\lambda \in S$. Then there are $a, b\in R^+$ such that $a<b$ and $[a,b]\subseteq S$.
   \end{lemma}

   \noindent{\bf Proof of Lemma} 
   We show
   \begin{equation}
     \lambda,\mu\in S \Rightarrow \lambda\mu\in S.
   \end{equation}
   Suppose $\lambda,\mu\in S$. In the following expression
   \begin{equation}\label{lmu}
     (f(\lambda\mu x) - f(x))\ln(x) = (f(\lambda\mu x) - f(\lambda x))\ln(x) + (f(\lambda x) - f(x))\ln(x)
   \end{equation}
   we have $\lim_{x\to \infty}(f(\lambda x) - f(x))\ln(x)=0$. Further,
   \begin{equation}
     \lim_{x\to \infty}(f(\lambda\mu x) - f(\lambda x))\ln(\lambda x)=\lim_{t\to \infty}(f(\mu t) - f(t))\ln(t)=0,
   \end{equation}
   and for $\lambda>1$, $\ln(x)< \ln(\lambda x)$, so
   \begin{equation}
     |(f(\lambda\mu x) - f(\lambda x))\ln(x)|\leq |(f(\lambda\mu x) - f(\lambda x))\ln(\lambda x)|,
   \end{equation}
   so $\lim_{x\to \infty}|(f(\lambda\mu x) - f(\lambda x))\ln(x)|=0$. By (\ref{lmu}) it follows
   \begin{equation}
     \lim_{x\to\infty}(f(\lambda\mu x) - f(x))\ln(x) = 0
   \end{equation}
   so, $\lambda\mu\in S$. Hence, $S$ is closed under multiplication, therefore, by Steinhaus lemma
   \cite{Steinhaus}, there are $a, b\in R^+$, $a<b$, so that $[a,b] \subseteq S$.
     \hfill $\Box$

  In the next lemma we show that the convergence interval $[a,b]$ can be expanded to $(0,\infty)$.
  \begin{lemma} \label{lemmaB}
    Let $f(x)$   be as in the theorem and suppose {\rm (\ref{llim})} for
    all $\lambda \in [a,b]$ for some $0<a<b$. Then {\rm (\ref{llim})} holds for all $\lambda\in R^+$.
  \end{lemma}

  \noindent{\bf Proof of Lemma}
  Let $\lambda\in[a,b]$ and $\mu\in R^+$ such that $a\leq \lambda/\mu \leq b$.
  Further,
  \begin{equation}\label{lmud}
     \begin{split}
        (f(\lambda x) - f(x))\ln(x) = &(f\left(\mu\frac{\lambda x}{\mu}\right) - f\left(\frac{\lambda x}{\mu}\right))\ln(x) + \\
                                      &(f\left(\frac{\lambda x}{\mu}\right) - f(x))\ln(x).
     \end{split}
  \end{equation}
  By assumptions on $\lambda$ and $\mu$, we have
  \begin{equation}
     \lim_{x\to \infty}(f(\lambda x) - f(x))\ln(x)=0, \quad \lim_{x\to \infty}(f\left(\frac{\lambda x}{\mu}\right) -
           f(x))\ln(x)=0,
  \end{equation}
  hence, by (\ref{lmud})
  \begin{equation}
    \lim_{x\to \infty}(f\left(\mu\frac{\lambda x}{\mu}\right) - f\left(\frac{\lambda x}{\mu}\right))\ln(x)=0
  \end{equation}
  For $x> \lambda/\mu$ we have
  \begin{equation}
     |(f\left(\mu\frac{\lambda x}{\mu}\right) - f\left(\frac{\lambda x}{\mu}\right))\ln(x)| \geq
     |(f\left(\mu\frac{\lambda x}{\mu}\right) - f\left(\frac{\lambda x}{\mu}\right))\ln(\lambda/\mu)|,
  \end{equation}
  so $\lim_{x\to \infty} (f\left(\mu\frac{\lambda x}{\mu}\right) - f\left(\frac{\lambda x}{\mu}\right))\ln(\lambda/\mu)=0$.
  But
  \begin{equation}
     \begin{split}
        \lim_{t\to \infty}(f(\mu t) - f(t))&\ln(t) =
            \lim_{x\to \infty}(f\left(\mu \frac{\lambda x}{\mu}\right) - f\left(\frac{\lambda x}{\mu}\right))
                     \ln\left(\frac{\lambda x}{\mu}\right) =  \\
        \lim_{x\to \infty}&(f\left(\mu \frac{\lambda x}{\mu}\right) - f\left(\frac{\lambda x}{\mu}\right))\ln(x)  + \\
             &\lim_{x\to \infty}(f\left(\mu \frac{\lambda x}{\mu}\right) - f\left(\frac{\lambda x}{\mu}\right))
                        \ln\left(\frac{\lambda}{\mu}\right)= 0
     \end{split}
  \end{equation}
  i.e.  $\lim_{t\to \infty}(f(\mu t) - f(t))\ln(t) =0$ for $a/b\leq \mu \leq b/a$. Hence we proved
  \begin{equation} \label{iterr}
    \begin{split}
      \bigwedge_{a\leq \lambda \leq b}\lim_{x\to \infty}(f(\lambda x) - f(x))\ln(x) =0 &\Rightarrow     \\
          \bigwedge_{a/b\leq \lambda \leq b/a}&\lim_{x\to \infty}(f(\lambda x) - f(x))\ln(x) =0
    \end{split}
  \end{equation}
  Iterating (\ref{iterr}) $n$ times, we obtain for arbitrary positive integer $n$
  \begin{equation}
      \bigwedge_{\left(\frac{a}{b}\right)^n\leq \lambda \leq \left(\frac{b}{a}\right)^n}\lim_{x\to \infty}(f(\lambda x) - f(x))\ln(x) =0.
  \end{equation}
  As $a/b<1$ and $b/a>1$ and so $\lim_{n\to\infty} (a/b)^n=0$ and $\lim_{n\to\infty} (b/a)^n=\infty$, we infer
  \begin{equation}
  \bigwedge_{\lambda \in R^+}\lim_{x\to \infty}(f(\lambda x) - f(x))\ln(x) =0.
  \end{equation}
    \hfill $\Box$

  Combining lemmas \ref{lemmaA} and \ref{lemmaB}, we obtain a proof of Theorem \ref{thmlm}.

  The following theorem gives us one interesting property of slowly varying functions.

  \begin{theorem}\label{SV}
  Assume $L(x)$ is a slowly varying function. Then there are measurable functions
  $\xi(x)\in \mathcal{Z}_0$ and $g(x)$, so that $L(x)=g(x)x^{\xi(x)}$,
  where  $g(x)\rightarrow g_0$ as $x \rightarrow \infty$, $g_0$ is a real positive constant.
  \end{theorem}

 \noindent{\bf Proof}. 
 Suppose that $L(x)$ is a slowly varying function.
  By integral representation theorem for SV functions, it follows that there are
  measurable functions $g(x)$, $\varepsilon\in \mathcal{Z}_0$ and $b\in \mathcal{R}$ so that
  \begin{equation}\label{RVrepresentation}
     L(x)= g(x) e^{\int_b^x\frac{\varepsilon(t)}{t}dt}, \quad x\geq b,
  \end{equation}
  and $ g(x)\to g_0$ as $x\to\infty$, $g_0$ is a real positive constant,
 wherefrom for $b=1$ we directly infer
 \begin{equation}
 L(x)=g(x)x^{\mathcal{L}(\varepsilon)(x)}.
 \end{equation}
 Since $\varepsilon(x)\in \mathcal{Z}_0$, by Proposition \ref{Llim}
 follows that $\mathcal{L}(\varepsilon)(x)\in \mathcal{Z}_0$, what proves the theorem.
  \hfill $\Box$
 \medskip

 The converse does not hold. If we suppose that $f(x)=g(x)x^{\xi(x)}$ is a SV function,
 where  $g(x)\rightarrow g_0$ as $x \rightarrow \infty$, $g_0$ is a real  positive constant, then for all $\lambda>0$ we have
 \begin{equation}\label{eq}
     1=\lim_{x\to\infty}{\dfrac{f(\lambda x)}{f(x)}}=
         \lim_{x\to\infty}{\dfrac{g(\lambda x)(\lambda x)^{\xi(\lambda x)}}{g(x)x^{\xi(x)}}}=
         \lim_{x\to\infty}{\dfrac{g(\lambda x){\lambda}^{\xi(\lambda x)} x^{\xi(\lambda x)}}{g(x)x^{\xi(x)}}}.
  \end{equation}
  Since $\lambda>0$, then $\lim_{x\to\infty} g(\lambda x)=\lim_{x\to\infty} g(x)=g_0$,
  and $\lim_{x\to\infty} \xi(\lambda x)=\lim_{x\to\infty} \xi(x)=0$.
  Therefore, the limit value in (\ref{eq}) depends only on the limit $\lim_{x\to\infty}{\dfrac{x^{\xi(\lambda x)}}{x^{\xi(x)}}}$.
  Furthermore we have
  \begin{equation}\label{ce}
    1= \lim_{x\to\infty}{\dfrac{x^{\xi(\lambda x)}}{x^{\xi(x)}}}=
       \lim_{x\to\infty} x^{\xi(\lambda x)-\xi(x)}=\lim_{x\to\infty} \exp((\xi(\lambda x)-\xi(x))\ln
       (x)),
  \end{equation}
  wherefrom  for every $\lambda>0$
  \begin{equation}\label{ce1}
    \lim_{x\to\infty}(\xi(\lambda x)-\xi(x))\ln (x)=0.
  \end{equation}
  However, that does not need to be the case, as the following example shows.
  Let $\xi(x)=\sin(x)/{\ln(x)}$ and $\lambda=\pi$. Taking the limit (\ref{ce1}) over  positive integers $n$, we obtain
  \begin{equation}
    \lim_{n\to\infty} \left( \sin(\pi n)/{\ln(\pi n)}-\sin(n)/{\ln(n)} \right) \ln (n)=
        \lim_{n\to\infty} \left( -\sin(n) \right),
  \end{equation}
  which does not exist, contradicting (\ref{ce1}).
  Therefore, not all functions \\
  $f(x)=g(x)x^{\xi(x)}$,
  where  $g(x)\rightarrow g_0$ as $x \rightarrow \infty$, $g_0$ is a real positive constant, have to be slowly varying.
  Next proposition gives a sufficient and necessary condition for representation  of normalized SV functions in Theorem \ref{SV}.

 \begin{theorem}\label{310} Let  $\xi\in \Z$. Then
    $ F(x)=x^{\xi(x)}$  is a SV function if and only if
    \begin{equation} \label{svv}
       \lim_{x\to\infty}(\xi(\lambda x)-\xi(x))\ln (x)=0
    \end{equation}
    uniformly with respect to $\lambda\in\R^+$.
  \end{theorem}

  \noindent{\bf Proof}. 
  First we prove Claim: Suppose $H(x,u)$ is a real function and \newline $\lim_{x\to \infty} H(x,u)= 1$
  uniformly with respect to $u\in I\subseteq \R$. Then
  \begin{equation}
     \lim_{x\to \infty} \ln(H(x,u))= 0
  \end{equation}
  uniformly with respect to $u\in I$.
  \smallskip

  \noindent Suppose $\lim_{x\to \infty} H(x,u)= 1$ uniformly with respect to $u\in I\subseteq \R$.
  Then for positive infinite x and $u\in I^*$, $H^*(x,u)= 1 + \varepsilon,\  \varepsilon \in \mu(0)$, and so
  $\ln(H^*(x,u))=   \varepsilon'$ for some $\varepsilon'\in \mu(0)$. Therefore, by Lemma \ref{unic}
  the convergence in question of $\ln(H(x,u))$ is uniform.
  Now we proceed to the proof of Theorem \ref{310}.
  \medskip

  \noindent $(\Rightarrow)$ Assume $ F(x)=x^{\xi(x)}$  is a SV function. Then
  \begin{equation}
     \lim_{x\to \infty} F(\lambda x)/F(x)= 1\quad  {\rm for}\hspace{2mm} {\rm all} \hspace{2mm} \lambda \in \R^+,
  \end{equation} i.e.
  $\lim_{x\to \infty} \lambda^{\xi(\lambda x)}e^{(\xi(\lambda x) -\xi(x))\ln(x)}=1$. By  Theorem \ref{karamata} (Karamata Uniform Convergence Theorem)
  it follows $\lim_{x\to \infty} e^{(\xi(\lambda x) -\xi(x))\ln(x)}=1$ uniformly with respect to $\lambda \in \R^+$.
  By Claim it follows that $\lim_{x\to \infty} (\xi(\lambda x) -\xi(x))\ln(x)= 0$ uniformly with respect to $\lambda \in \R^+$.
  \smallskip

  \noindent $(\Leftarrow)$  Assume (\ref{svv}) uniformly with respect to $\lambda\in\R^+$.  Then
  for positive infinite $x$ and $\lambda \in R^{+*}$, $\xi(\lambda x)= \eta$ and due to the uniform convergence
  by Lemma \ref{unic}, $(\xi(\lambda x) - \xi(x))\ln(x)=\varepsilon$ are infinitesimals, so
  \begin{equation}
     \frac{F(\lambda x)}{F(x)}= \lambda^{\xi(\lambda x)} e^{(\xi(\lambda x) - \xi(x))\ln(x)}=
           \lambda^\eta e^\varepsilon \approx 1.
  \end{equation}
  Hence $\lim_{x\to \infty} F(\lambda x)/F(x) = 1$ and in fact, by Theorem \ref{karamata}, this convergence is uniform
  with respect to $\lambda\in R^+$.
  \hfill $\Box$
  \smallskip

  With a simple modification of above proof, one can prove a variant of the previous theorem for
  functions of the form $F(x)= g(x) x^{\xi(x)}$, where $g(x) \to 1$ as $x\to \infty$ and $\xi\in \Z$.

\section*{Conclusion}

We used methods of nonstandard analysis in order to obtain proof of a generalization of Karamata uniform convergence theorem for slowly varying functions. We introduced operator $\mathcal{L}$ and proved its several properties. Furthermore, the connection between the operator $\mathcal{L}$ and slowly varying functions is derived. Moreover, some properties of slowly varying functions are obtained.

\bibliographystyle{plain}
\bibliography{ArXiv_A_nonstandard_approach_to_Karamata_uniform_convergence_theorem}

\begin{thebibliography}{10}

\bibitem{bingham}
N.H. Bingham, C.M. Goldie, and J.L. Teugels.
\newblock {\em Regular Variation}.
\newblock Encyclopedia of Mathematics and its Applications. Cambridge
  University Press, 1987.

\bibitem{Keisler}
C.C. Chang and H.J. Keisler.
\newblock {\em Model Theory}.
\newblock North Holland, 3rd edition, 1990.

\bibitem{Davis}
M.~Davis.
\newblock {\em Applied {N}onstandard {A}nalysis}.
\newblock Dover Publications, 2005.

\bibitem{karamata}
J.~Karamata.
\newblock Sur une mode de croissance r\'eguliere des functions.
\newblock {\em Math.(Cluj)}, 4:38--53, 1930.

\bibitem{liddle}
A.R. Liddle and D.H. Lyth.
\newblock {\em Cosmological Inflation and Large-Scale Structure}.
\newblock Cambridge Univ. Press, 2000.

\bibitem{Seneta}
E.~Seneta.
\newblock {\em Regularly Varying Functions}.
\newblock Springer Berlin, Heidelberg, 1 edition, 1976.

\bibitem{Steinhaus}
Hugo Steinhaus.
\newblock Sur les distances des points dans les ensembles de mesure positive.
\newblock {\em Fundamenta Mathematicae}, 1(1):93--104, 1920.

\bibitem{stern}
I.~{Stern}.
\newblock {On Fractal Modeling in Astrophysics: The Effect of Lacunarity on the
  Convergence of Algorithms for Scaling Exponents.}
\newblock In Gareth {Hunt} and Harry {Payne}, editors, {\em Astronomical Data
  Analysis Software and Systems VI}, volume 125 of {\em Astronomical Society of
  the Pacific Conference Series}, pages 222--225, January 1997.

\bibitem{Stroyan}
K.D. Stroyan and W.A.J. Luxemburg.
\newblock {\em Introduction to the Theory of Infinitesimals}.
\newblock Academic Press, New York, 1976.

\bibitem{mijajlo2014}
\v{Z}. Mijajlovi\'{c}, D.~Aran{\dj}elovi\'{c}, M.~Ra\v{s}kovi\'{c}, and
  R.~{\Dj}or{\dj}evi\'{c}.
\newblock {\em {Nestandardna Analiza}}.
\newblock {Matemati\v{c}ki fakultet Beograd}, 2014.

\bibitem{mijajlo2007}
\v{Z}. Mijajlovi\'{c} and N.~Pejovi\'{c}.
\newblock Non‐archimedean methods in cosmology.
\newblock {\em AIP Conference Proceedings}, 895(1):317--320, 2007.

\bibitem{mijajlo2015}
\v{Z}. Mijajlovi\'{c}, N.~Pejovi\'{c}, and V.~Mari\'{c}.
\newblock {On the $\varepsilon$ Cosmological Parameter}.
\newblock {\em Serbian Astronomical Journal}, 190:25--31, June 2015.

\bibitem{mijajlo2007A}
\v{Z}. Mijajlovi\'{c}, N.~Pejovi\'{c}, and S.~Ninkovi\'{c}.
\newblock {Nonstandard Representations of Processes in Dynamical Systems}.
\newblock {\em AIP Conference Proceedings}, 934(1):151--157, 2007.

\bibitem{mijajlo2019}
\v{Z}. Mijajlovi\'{c}, N.~Pejovi\'{c}, and V.~Radovi\'{c}.
\newblock Asymptotic solution for expanding universe with matter-dominated
  evolution.
\newblock {\em International Journal of Geometric Methods in Modern Physics},
  16(04):1950063, 2019.

\bibitem{mijajlo2012}
\v{Z}. Mijajlovi\'{c}, N.~Pejovi\'{c}, S.~\v{S}egan, and G.~Damljanovi\'{c}.
\newblock On asymptotic solutions of {F}riedmann equations.
\newblock {\em Applied Mathematics and Computation}, 219(3):1273--1286, 2012.

\end{thebibliography}

  \end{document}